\documentclass[reqno,a4paper]{article}
\usepackage[english]{babel}
\usepackage{amsmath}
\usepackage{amsthm}
\usepackage{amssymb}
\usepackage{enumerate}
\usepackage{xspace}
\usepackage{euscript}
\usepackage{graphicx}
\usepackage{graphics}
\usepackage{amscd}
\usepackage{epsfig}
\usepackage{tabularx}
\usepackage[usenames]{color}
\newtheorem{lem}{Lemma}[section]
\newtheorem{defn}{Definition}[section]
\newtheorem{thm}[lem]{Theorem}
\newtheorem{rem}[lem]{Remark}

\newtheorem{asm}[lem]{Assumption}
\numberwithin{equation}{section}
\newcommand{\mri}{{\mathrm{i}}} 

\title{\textbf{Lipschitz stability at the boundary for time-harmonic diffuse optical tomography}}
\author{\textsc{Olga Doeva\thanks{Bernal Institute,
University of Limerick, Castletroy, Limerick, Ireland,
Olga.Doeva$@$ul.ie},\quad Romina
Gaburro}\thanks{Department of Mathematics and Statistics, Health Research Institute (HRI),
University of Limerick, Ireland,
Romina.Gaburro$@$ul.ie},\\ \textsc{William R. B. Lionheart\thanks{Department of Mathematics, University of Manchester, Alan Turing Building,Oxford Rd Manchester M13 9PL, United Kingdom, bill.lionheart@manchester.ac.uk},\quad and \quad Clifford J. Nolan\thanks{Department of Mathematics and Statistics,
Health Research Institute (HRI), University of Limerick, Ireland, Clifford.Nolan$@$ul.ie}}}

\date{}

\begin{document}
\maketitle

\par{\centering \textit{This manuscript is dedicated to Sergio Vessella on the occasion of his $65^{th}$ birthday, to honour his outstanding contribution to the fields of inverse problems and mathematical analysis.}\par}

 \begin{abstract}
We study the inverse problem in Optical Tomography of determining the optical properties of a medium $\Omega\subset\mathbb{R}^n$, with $n\geq 3$, under the so-called \textit{diffusion approximation}. We consider the time-harmonic case where $\Omega$ is probed with an input field that is modulated with a fixed harmonic frequency $\omega=\frac{k}{c}$, where $c$ is the speed of light and $k$ is the wave number. We prove a result of Lipschitz stability of the \textit{absorption coefficient} $\mu_a$ at the boundary $\partial\Omega$ in terms of the measurements in the case when the \textit{scattering coefficient} $\mu_s$ is assumed to be known and $k$ belongs to certain intervals depending on some \textit{a-priori} bounds on $\mu_a$, $\mu_s$.
\end{abstract}


\section{Introduction}\label{section introduction}
\setcounter{equation}{0}

Although Maxwell's equations provide a complete model for the light propagation in a scattering medium on a micro scale, on the scale suitable for medical diffuse Optical Tomography (OT) an appropriate model is given by the \textit{radiative transfer equation} (or \textit{Boltzmann equation})\cite{ArSc}. If $\Omega$ is a domain in $\mathbb{R}^{n}$, with $n\geq 2$ with smooth boundary $\partial\Omega$ and radiation is considered in the body $\Omega$, then it is well known that if the input field is modulated with a fixed harmonic frequency $\omega$, the so-called \textit{diffusion approximation} leads to the complex partial differential equation (see \cite{Ar}) for the energy current density $u$

\begin{equation}\label{K,q equation}
-\mbox{div}\left(K\nabla u\right)+(\mu_a -\mri k)u=0,\qquad\textnormal{in}\quad\Omega.
\end{equation}

Here $k=\frac{\omega}{c}$ is the wave number, $c$ is the speed of light and, in the anisotropic case, the so-called \textit{diffusion tensor} $K$, is the complex matrix-valued function

\begin{equation}\label{K time-harmonic 1}
K=\frac{1}{n}\Big((\mu_a-\mri k) I+(I-B)\mu_s\Big)^{-1},\qquad\textnormal{in}\quad\Omega,
\end{equation}

where $B_{ij}=B_{ji}$ is a real matrix-valued function, $I$ is the $n\times n$ identity matrix and $I-B$ is positive definite (\cite{Ar}, \cite{HS}, \cite{H}) on $\Omega$. The spacially dependent real-valued coefficients $\mu_a$ and $\mu_s$ are called the \textit{absorption} and the \textit{scattering coefficients} of the medium $\Omega$ respectively and represent the optical properties of $\Omega$. It is worth noticing that many tissues including parts of the brain, muscle and breast tissue have fibrous structure on a microscopic scale which results in anisotropic physical properties on a larger scale. Therefore the model considered in this manuscript seems appropriate for the case of medical applications of OT (see \cite{HAS}). Although it is common practise in OT to use the Robin-to-Robin map to describe the boundary measurements (see \cite{Ar}), the Dirichlet-to-Neumann (D-N) map will be employed here instead. This is justified by the fact that in OT, prescribing its inverse, the Neumnann-to-Dirichlet (N-D) map (on the appropriate spaces), is equivalent to prescribing the Robin-to-Robin boundary map. A rigorous definition of the D-N map for equation \eqref{K,q equation} will be given in section \ref{formulation problem}.



It is also well known that prescribing the N-D map is insufficient to recover
both coefficients $\mu_a$ and $\mu_s$ uniquely \cite{ArL} unless \textit{a-priori} smoothness assumptions are imposed \cite{Ha}. In this paper we consider the problem of determining the absorption coefficient $\mu_a$ in a medium  $\Omega\subset\mathbb{R}^n$, $n\geq 3$, that is probed with an input field which is modulated with a fixed harmonic frequency $\omega=\frac{k}{c}$, with $k\neq 0$ (time-harmonic case) and whose scattering coefficient $\mu_s$ is assumed to be known. More precisely, we show that $\mu_a$, restricted to the boundary $\partial\Omega$, depends upon the D-N map of \eqref{K,q equation}, $\Lambda_{K, \mu_a}$, in a Lipschitz way when $k$ is chosen in certain intervals that depend on \textit{a-priori} bounds on $\mu_a$, $\mu_s$ and on the ellipticity constant of $I-B$ (Theorem  \ref{stabilita' al bordo}). The static case ($k=0$), for which \eqref{K,q equation} is a single real elliptic equation, was studied in \cite{G}, where the author proved Lipschitz stability of $\mu_a$ and H\"older stability of the derivatives  of $\mu_a$ at the boundary in terms of $\Lambda_{K,\mu_a}$. In the present paper we show that in the time-harmonic case, for which \eqref{K,q equation} is a complex elliptic equation, a Lipschitz stability estimate of $\mu_a$ at the boundary $\partial\Omega$ in terms of $\Lambda_{K,\mu_a}$ still holds true if $k$ is chosen within certain ranges. The case where $\mu_a$ is assumed to be known and the scattering coefficient $\mu_s$ is to be determined, can be treated in a similar manner. The choice  in this paper of focusing on the determination of $\mu_a$ rather than the one of $\mu_s$ is driven by the medical application of OT we have in mind. While $\mu_s$ varies from tissue to tissue, it is the absorption coefficient $\mu_a$ that carries the more interesting physiological information as it is related to the global concentrations of certain metabolites in their oxygenated and deoxygenated states.

Our main result (Theorem \ref{stabilita' al bordo}) is based on the construction of singular solutions to the complex elliptic equation \eqref{K,q equation}, having an isolated singularity outside $\Omega$. Such solutions were first constructed in \cite{A1} for equations of the type

\begin{equation}\label{Calderon operator}
\mbox{div}(K\nabla u ) = 0,\qquad\textnormal{in}\quad\Omega,
\end{equation}

when $K$ is a real matrix-valued function belonging to $W^{1,p}(\Omega)$, with $p>n$ and they were employed to prove stability results at the boundary in \cite{A1}, \cite{AG}, \cite{AG1} and \cite{GL} in the case of Calder\'on's problem (see \cite{C}) with global, local data and on manifolds. The singular solutions introduced in \cite{A1} were extended in \cite{Sa} to equations of the type

\begin{equation}\label{Salo operator}
-\mbox{div}(K\nabla u + Pu) + Q\cdot\nabla u +qu =0,\qquad\textnormal{in}\quad\Omega,
\end{equation}

with real coefficients, where $K$ is merely H\"older continuous. Singular solutions were also studied in \cite{I1}.

In this paper we extend the singular solutions introduced in \cite{A1} to the case of elliptic equations of type \eqref{K,q equation} with complex coefficients. Such a construction is done by treating \eqref{K,q equation} as a strongly elliptic system with real coefficients, since $\Re K\geq\tilde\lambda^{-1}I>0$, where  $\tilde\lambda$ is a positive constant depending on the \textit{a-priori} information on $\mu_s$, $B$ and $\mu_a$. We wish to stress, however, that in \cite{A1} the author constructed singular solutions to \eqref{Calderon operator} which have an isolated singularity of arbitrary high order, where the current paper extends such construction to singular solutions to the complex equation \eqref{K,q equation} having an isolated singularity of Green's type only. This is sufficient to prove the Lipschitz continuity of the boundary values of $\mu_a$ in terms of the D-N map. The more general construction of the singular solutions with an isolated singularity of arbitrary high order for elliptic complex partial differential equations will be material of future work. 

This paper is stimulated by the work of Alessandrini and Vessella \cite{A-V}, where the authors proved global Lipschitz stability of the conductivity in a medium $\Omega$ in terms of the D-N map for Calder\'on's problem, in the case when the conductivity is real, isotropic and piecewise constant on a given partition of $\Omega$. This fundamental result was extended to the complex case in \cite{Be-Fr} and in the context of various inverse problems for example in \cite{A-dH-G-S}, \cite{A-dH-G-S1}, \cite{Be-dH-Fr-V-Z} and  \cite{A-dH-G}, \cite{A-dH-G-S2}, \cite{G-S} in the isotropic and anisotropic settings respectively. The machinery of the proof introduced in \cite{A-V} is based on an induction argument that combines quantitative estimates of unique continuation together with a careful asymptotic analysis of Green's functions. The initial step of their induction argument relies on Lipschitz (or H\"older) stability estimates at the boundary of the physical parameter that one wants to estimate in terms of the boundary measurements, which is the subject of the current manuscript. Our paper also provides a first step towards a reconstruction procedure of $\mu_a$ by boundary measurements based on a Landweber iterative method for non-linear problems studied in \cite{dHQS}, where the authors provided an analysis of the convergence of such algorithm in terms of either a H\"older or Lipschitz global stability estimates (see also \cite{A-dH-F-G-S}). We also refer to \cite{KVKaAr} and \cite{HebdenArridge} for further reconstruction techniques of the optical properties of a medium.

The paper is organized as follows. Section 2 contains the formulation of the problem (subsections \ref{main assumptions}, \ref{D-to-N}) and our main result (subsection \ref{sub main result}, Theorem \ref{stabilita' al bordo}). Section 3 is devoted to the construction of singular solutions of equation \eqref{K,q equation} having a Green's type isolated singularity outside $\Omega$. The proof of our main result (Theorem \ref{stabilita' al bordo}) is given in section \ref{proofs main result}.



\section{Formulation of the problem and main result}\label{formulation problem}

\subsection{Main assumptions}\label{main assumptions}
We rigorously formulate the problem by introducing the following notation, definitions and assumptions.
For $n\geq 3$, a point $x\in \mathbb{R}^n$ will be denoted by $x=(x',x_n)$, where $x'\in\mathbb{R}^{n-1}$ and $x_n\in\mathbb{R}$.
Moreover, given a point $x\in \mathbb{R}^n$, we will denote with $B_r(x), B_r'(x')$ the open balls in
$\mathbb{R}^{n},\mathbb{R}^{n-1}$, centred at $x$ and $x'$ respectively with radius $r$
and by $Q_r(x)$ the cylinder

\[Q_r(x)=B_r'(x')\times(x_n-r,x_n+r).\]

We will also denote $B_r=B_r(0)$, $B'_r=B'_r(0)$ and $Q_r=Q_r(0)$.



\begin{defn}\label{def boundary}
Let $\Omega$ be a bounded domain in $\mathbb R^n$, with $n\geq 3$. We shall say that the boundary of $\Omega$, $\partial\Omega$, is of Lipschitz class with constants $r_0,L>0$, if for any $P\in\partial\Omega$ there exists a rigid
transformation of coordinates under which we have $P=0$ and

$$\Omega\cap Q_{r_0}=\{(x',x_n)\in Q_{r_0}\: |\,x_n>\varphi(x')\},$$

where $\varphi$ is a Lipschitz function on $B'_{r_0}$ satisfying

\[\varphi(0)=0\]

and

\[\|\varphi\|_{C^{0,1}(B'_{r_0})}\leq Lr_0.\]

\end{defn}

We consider, for a fixed $k>0$,

\begin{equation}\label{L complex}
L= - \mbox{div}\left(K\nabla\cdot\right) + q,\qquad\textnormal{in}\quad\Omega,
\end{equation}

where $K$ is the complex matrix-valued function

\begin{equation}\label{K time-harmonic 2}
K(x)=\frac{1}{n}\Big((\mu_a(x)-\mri k) I+(I-B(x))\mu_s(x)\Big)^{-1},\qquad\textnormal{for\:any}\:x\in\Omega,
\end{equation}

and $q$  is the complex-valued function 

\begin{equation}\label{q}
q=\mu_a - ik\qquad\textnormal{in}\quad\Omega.
\end{equation}

We recall that $I$ denotes the $n\times n$ identity matrix, where the matrix $B$ is given by the OT physical experiment and it is such that $B\in L^{\infty}(\Omega, Sym_n)$, where $Sym_n$ denotes the class of $n\times n$ real-valued symmetric matrices and such that $I-B$ is a positive definite matrix (\cite{Ar}, \cite{HS}, \cite{H}). In this paper we assume that the \textit{scattering coefficient} $\mu_s$ is also known in $\Omega$ and it is the \textit{absorption coefficient} $\mu_a$ that we seek to estimate from boundary measurements.\\

We assume that there are positive constants $\lambda$, $E$ and $\mathcal{E}$ and $p>n$ such that the known quantities $B$, $\mu_s$ and the unknown quantity $\mu_a$ satisfy the two assumptions below respectively.

\begin{asm}\textnormal{(Assumption on $\mu_s$ and $B$)}\label{assumption on mus and B}\\

\begin{equation}\label{assumption on scattering coeff positive}
\lambda^{-1}\leq\mu_{s}(x)\leq \lambda,\qquad\textnormal{for\:a.e.}\: x\in\Omega,
\end{equation}

\begin{equation}\label{assumption scattering sobolev}
||\mu_s||_{W^{1,\:p}(\Omega)} \leq E
\end{equation}

and 

\begin{equation}
\mathcal{E}^{-1}|\xi|^2\leq(I-B(x))\xi\cdot\xi\leq\mathcal{E}|\xi|^2,\qquad\textnormal{for\:a.e.\:}x\in\Omega,\quad\textnormal{for\:any}\:\xi\in\mathbb{R}^n.
\end{equation}

\end{asm}

\begin{asm}\textnormal{(Assumption on $\mu_a$)}\label{assumption on mua}\\

\begin{equation}\label{limitazioni per a,b}
\lambda^{-1}\leq\mu_{a}(x),\leq\lambda,\qquad\textnormal{for\:a.e.}\: x\in\Omega,
\end{equation}

\begin{equation}\label{holder a,b}
\parallel\mu_{a}\parallel_{\:W^{1,p}(\Omega)}\leq{E}.
\end{equation}



\end{asm}

We state below some facts needed in the sequel of the paper. Most of them are straightforward consequences of our assumptions.\\

The inverse of $K$

\begin{equation}\label{K inverse}
K^{-1}= n \Big(\mu_a I + (I-B)\mu_s -ikI\Big),\qquad\textnormal{on}\quad\Omega
\end{equation}

has real and imaginary parts given by the symmetric, real matrix valued-functions on $\Omega$

\begin{eqnarray}
K^{-1}_{R} &=& n\left(\mu_a I + (I-B)\mu_s\right),\label{K inverse 1}\\
K^{-1}_I &=& \!\!\!\!-nkI\label{K inverse 2}
\end{eqnarray}

respectively. As an immediate consequence of assumptions \ref{assumption on mus and B}, \ref{assumption on mua} we have

\begin{eqnarray}\label{apriori assumptions inequalities}
n\lambda^{-1}(1+\mathcal{E}^{-1})|\xi|^2\leq K^{-1}_R (x) \xi\cdot\xi &\leq & n\lambda(1+\mathcal{E})|\xi|^2,\label{apriori ass ineq 1}\\
- K^{-1}_I (x)\xi \cdot \xi &=& nk|\xi|^2,\label{apriori ass ineq 2}
\end{eqnarray}

for a.e. $x\in\Omega$ and any $\xi\in\mathbb{R}^n$. Moreover $K^{-1}_{R}$ and $K^{-1}_I$  commute, therefore the real and imaginary parts of $K$ are the symmetric, real matrix valued-functions on $\Omega$

\begin{eqnarray}
K_R &=& \frac{1}{n}\bigg(\Big( \mu_a I + (I-B)\mu_s \Big)^2 + k^2I\bigg)^{-1}\big(\mu_a I + (I-B)\mu_s \big)\label{K1},\\
K_I  &=& \frac{k}{n}\bigg(\Big( \mu_a I + (I-B)\mu_s \Big)^2 + k^2I\bigg)^{-1}\label{K2}
\end{eqnarray}

respectively. Assumptions \ref{assumption on mus and B}, \ref{assumption on mua} also imply that  

\begin{eqnarray}
& & K_R(x) \xi\cdot\xi \geq\frac{\lambda (1+\mathcal{E})}{n} \Big(\lambda^2 (1+\mathcal{E})^2 + k^2\Big)^{-1} |\xi|^2 ,\label{K1 pos def}\\
& & K_I(x) \xi\cdot\xi \geq\frac{k}{n} \Big(\lambda^{2} (1+\mathcal{E})^2 + k^2 \Big)^{-1} |\xi|^2 \label{K2 neg def},
\end{eqnarray}

for a.e. $x\in\Omega$, for every $\xi\in\mathbb{R}^n$ and the \textit{boundness condition}

\begin{equation}\label{boundness OT}
|K_R(x)|^2 + |K_I(x)|^2\leq \Big(\lambda^{-2} (1+\mathcal{E}^{-1})^2 + k^2\Big)^{-2}\:\Big(\frac{\lambda^2 (1+\mathcal{E})^2 + k^2}{n^2}\Big), 
\end{equation}

for a.e. $x\in\Omega$.

Moreover $K=\{K^{hk}\}_{h,k=1,\dots ,n}$ and $q$ satisfy

\begin{equation}\label{boundness K}
||K^{hk}||_{W^{1,p}(\Omega)}\leq C_1,\qquad h,k=1,\dots , n,
\end{equation}

and 

\begin{equation}\label{boundness q}
|q(x)| = |\mu_a (x) -ik|\leq \lambda +k,\qquad\textnormal{for\:a.e.}\:x\in\Omega,
\end{equation}

respectively, where $C_1$ is a positive constant depending on $\lambda$, $E$, $\mathcal{E}$, $k$ and $n$.


By denoting $q=q_R+iq_I$, the complex equation

\begin{equation}\label{complex eq}
-\mbox{div}\left(K\nabla u\right)+qu=0,\qquad\textnormal{in}\:\Omega
\end{equation}

is equivalent to the system for the vector field $u=(u^1, u^2)$

\begin{equation}\label{system OT}
\left\{ \begin{array}{ll} - \textnormal{div}(K_R\nabla
u^1) +  \textnormal{div}(K_I\nabla
u^2) + \left(q_Ru^1-q_Iu^2\right)=0, &
\textrm{$\textnormal{in}\quad\Omega$},\\
-  \textnormal{div}(K_I\nabla
u^1) +  \textnormal{div}(K_R\nabla
u^2 ) + \left(q_Iu^1+q_R u^2\right)=0,&
\textrm{$\textnormal{in}\quad\Omega$},
\end{array} \right.
\end{equation}

which can be written in a more compact form as

\begin{equation}\label{system compact}
-\mbox{div}(C\nabla u) + qu = 0,\qquad\textnormal{in}\quad\Omega
\end{equation}

or, in components, as

\begin{equation}\label{system compact components}
-\frac{\partial}{\partial x_h}\left\{C_{lj}^{hk} \frac{\partial}{\partial x_k} u^j \right\} + q_{lj}u^j = 0,\qquad\textnormal{for}\quad l=1,2, \quad\textnormal{in}\quad\Omega,
\end{equation}

where $\left\{C_{lj}^{hk}\right\}_{h,k=1,\dots, n}$ is defined by

\begin{equation}\label{def C}
C_{lj}^{hk}=K_R^{hk}\delta_{lj}-K_I^{hk}\left(\delta_{l1}\delta_{j2} - \delta_{l2}\delta_{j1}\right)
\end{equation}

and $\{q_{lj}\}_{l,j=1,2}$ is a $2\times 2$ real matrix valued function on $\Omega$ defined by

\begin{equation}\label{def q}
q_{lj}=q_R\delta_{lj}-q_I\left(\delta_{l1}\delta_{j2} - \delta_{l2}\delta_{j1}\right).
\end{equation}

\eqref{K1 pos def}, together with \eqref{boundness OT} imply that system \eqref{system OT} is \textit{uniformly elliptic} and \textit{bounded}, therefore it satisfies the \textit{strong ellipticity condition}

\begin{equation}\label{strong ellipticity}
C_2^{-1}|\xi|^2\leq C_{lj}^{hk}(x)\xi_{h}^{l}\xi_{k}^{j}\leq C_2 |\xi|^2, \qquad\textnormal{for\:a.e.}\: x\in\Omega ,\quad \textnormal{for\:all}\:\xi\in\mathbb{R}^{2n},
\end{equation}

where $C_2 >0$ is a constant depending on $\lambda$, $\mathcal{E}$, $k$ and $n$.

\begin{rem}
Matrix  $q =\{q_{lj}\}_{l,j=1}^{2}$ 

\begin{equation}
\begin{pmatrix}
\:\mu_a & k \\
\!\!\!- k & \mu_a
\end{pmatrix}
\end{equation}

is uniformly positive definite on $\Omega$ and it satisfies

\begin{equation}\label{q positive definite}
\lambda^{-1}|\xi|^2\leq q(x)\xi\cdot\xi\leq\lambda |\xi|^2,\qquad\textnormal{for\:a.e.}\: x\in\Omega,\quad\textnormal{for\:every}\:\xi\in\mathbb{R}^2.
\end{equation}

\end{rem}

\begin{defn}\label{a priori data}
We will refer in the sequel to the set of positive numbers $r_0$, $L$, $\lambda$, $E$, $\mathcal{E}$ introduced above, along with
the space dimension $n$, $p>n$, the wave number $k$ and the diameter of $\Omega$, $diam(\Omega)$, as to the \textit{a-priori data}. 
\end{defn}


\subsection{The Dirichlet-to-Neumann map.}\label{D-to-N}
Let $K$ be the complex matrix valued - function on $\Omega$ introduced in \eqref{K time-harmonic 2} and $q=\mu_a - ik$,  satisfying assumptions \ref{assumption on mus and B}, \ref{assumption on mua}. $B$ and $\mu_s$ are assumed to be known in $\Omega$ and satisfying assumption \ref{assumption on mus and B}, so that $K$ is completely determined by $\mu_a$, satisfying assumption \ref{assumption on mua}, on $\Omega$. Denoting by $\langle\cdot,\cdot\rangle$ the $L^{2}(\partial\Omega)$-pairing between $H^{\frac{1}{2}}(\partial\Omega)$ and its dual $H^{-\frac{1}{2}}(\partial\Omega)$, we will emphasise such dependence of $K$ on $\mu_a$ by denoting $K$ by
 
\[K_{\mu_a}.\]

For any $v,w\in\mathbb{C}^n$, with $v=(v_1,\dots , v_n)$, $w=(w_1,\dots , w_n)$, we will denote throughout this paper by $v\cdot w$, the expression
\[v\cdot w = \sum_{i=1}^{n} v_i w_i.\]

\begin{defn}
The Dirichlet-to-Neumann (D-N) map corresponding to $\mu_a$ is the operator

\begin{equation}\label{mappaDN}
\Lambda_{\mu_a}:H^{\frac{1}{2}}(\partial\Omega)\longrightarrow{H}^{-\frac{1}{2}}(\partial\Omega)
\end{equation}

defined by

\begin{equation}\label{def DN}
\langle\Lambda_{\mu_a}\:f,\:\overline{g}\rangle\:=\:\int_{\:\Omega}\Big( K_{\mu_a}(x) \nabla{u}(x)\cdot\nabla\varphi(x)+(\mu_a(x)-ik)u(x)\varphi(x)\Big)\:dx,
\end{equation}

for any $f$, $g\in H^{\frac{1}{2}}(\partial\Omega)$, where $u\in{H}^{1}(\Omega)$ is the weak solution to

\begin{displaymath}
\left\{ \begin{array}{ll} -\textnormal{div}(K_{\mu_a}(x)\nabla
u(x))+ (\mu_a-ik)(x)u(x)=0, &
\textrm{$\textnormal{in}\quad\Omega$},\\
u=f, & \textrm{$\textnormal{on}\quad{\partial\Omega}$}
\end{array} \right.
\end{displaymath}

and $\varphi\in H^{1}(\Omega)$ is any function such that $\varphi\vert_{\partial\Omega}=g$ in the trace sense.
\end{defn}


Given $B$, $\mu_s$, $\mu_{a_i}$, and the corresponding diffusion tensors $K_{\mu_{a_i}}$, for $i=1,2$, satisfying assumptions \ref{assumption on mus and B}, \ref{assumption on mua}, the well known Alessandrini's identity (see \cite[(5.0.4), p.129]{A1})

\begin{eqnarray}\label{Alessandrini identity}
\langle \left(\Lambda_{\mu_{a_1}} - \Lambda_{\mu_{a_2}}\right) f,\overline{g}\rangle &=& \int_{\Omega} \left(K_{\mu_{a_1}}(x) - K_{\mu_{a_2}}(x)\right)\nabla u(x)\cdot\nabla v(x)\:dx\nonumber\\
&+&\int_{\Omega}\left(\mu_{a_1}(x)-\mu_{a_2}(x)\right)u(x)v(x)\:dx,
\end{eqnarray}

holds true for any $f, g\in H^{\frac{1}{2}}(\partial\Omega)$, where $u, v\in H^{1}(\Omega)$ are the unique weak solutions to the Dirichlet problems

\begin{displaymath}
\left\{ \begin{array}{ll} -\textnormal{div}(K_{\mu_{a_1}}(x)\nabla
u(x)) +  (\mu_{a_1}-ik)u(x)=0, &
\textrm{$\textnormal{in}\quad\Omega$},\\
u=f, & \textrm{$\textnormal{on}\quad{\partial\Omega}$}
\end{array} \right.
\end{displaymath}
and

\begin{displaymath}
\left\{ \begin{array}{ll} -\textnormal{div}(K_{\mu_{a_2}}(x)\nabla
v(x)) +  (\mu_{a_2}-ik)v(x)=0, &
\textrm{$\textnormal{in}\quad\Omega$},\\
v=g, & \textrm{$\textnormal{on}\quad{\partial\Omega}$},
\end{array} \right.
\end{displaymath}

respectively.\\

We will denote in the sequel by $\parallel\cdot\parallel_{\mathcal{L}(H^{\frac{1}{2}}(\partial\Omega),H^{-\frac{1}{2}}(\partial\Omega))}$ the norm on the Banach space of
bounded linear operators between $H^{\frac{1}{2}}(\partial\Omega)$ and $H^{-\frac{1}{2}}(\partial\Omega)$.


\subsection{The main result}\label{sub main result}

%
%

\begin{thm}\label{stabilita' al bordo}(\textnormal{Lipschitz stability of boundary
values}). Let $n\geq 3$, and $\Omega$ be a bounded domain in $\mathbb{R}^n$ with
Lipschitz boundary with constants $L, r_0$ as in definition \ref{def boundary}. If $p>n$, $B$, $\mu_{s}$ and $\mu_{a_{i}}$, for $i=1,2$, satisfy assumptions \ref{assumption on mus and B}, \ref{assumption on mua} and the wave number $k$ satisfies either

\begin{equation}\label{k range 1}
0<k \leq k_0 :=\frac{ \sqrt{\lambda^2 (1+\mathcal{E})^2 + \lambda^{-2} (1+\mathcal{E}^{-1})^{2} \tan^{2}\left(\frac{\pi}{2n}\right)} -\lambda (1+\mathcal{E})}{\tan\left(\frac{\pi}{2n}\right)},
\end{equation}

or

\begin{equation}\label{k range 2}
k \geq \tilde{k}_0 :=\frac{1 +  \sqrt{1 + \tan^{2}\left(\frac{\pi}{2n}\right)}}{\tan\left(\frac{\pi}{2n}\right)}\: \lambda (1+\mathcal{E}),
\end{equation}

where, $\lambda$ and $\mathcal{E}$ are the positive numbers introduced in assumptions \ref{assumption on mus and B}, \ref{assumption on mua}, then 

\begin{equation}\label{stabilita' anisotropa}
\parallel\mu_{a_{1}}-\mu_{a_{2}}\parallel_{L^{\infty}\:(\partial\Omega)}
\leq{C}\parallel\Lambda_{\mu_{a_1}}-\Lambda_{{\mu_{a_2}}}\parallel_{\mathcal{L}(H^{\frac{1}{2}}(\partial\Omega),H^{-\frac{1}{2}}(\partial\Omega))},
\end{equation}

\noindent where $C>0$ is a constant depending on $n$, $p$, $L$, $r_0$, $diam(\Omega)$, $\lambda$, $E$, $\mathcal{E}$ and $k$.
\end{thm}

\section{Singular solutions}\label{section singular solutions bis}

We consider 

\begin{equation}\label{L complex 2}
L= - \mbox{div}\left(K\nabla\cdot\right) +q,\qquad\textnormal{in}\quad B_R=\Big\{x\in\mathbb{R}^n\:  \big | \: |x|<R\Big\},
\end{equation}

where $K=\{K^{hk}\}_{h,k=1,\dots , n}$ and $q$ are the complex matrix valued-function and the complex function respectively introduced in section \ref{section introduction} and satisfying assumptions \ref{assumption on mus and B}, \ref{assumption on mua} on $B_R$.

\begin{thm}\textnormal{(Singular solutions for $L=-\textnormal{div}( K\nabla\cdot)+q$)}.\label{theor singular sol}
Given $L$ on $B_R$ as in \eqref{L complex 2}, there exists $u\in{W}_{loc}^{2,\:p}(B_{R}\setminus\{0\})$ such that

\begin{equation}\label{solution with potential}
Lu=0,\mbox{ in } {B_R\setminus\{0\}}
\end{equation}

and furthermore

\begin{equation}\label{u singular}
u(x) = \left(K^{-1}(0)x\cdot x\right)^{\frac{2-n}{2}} + w(x),
\end{equation}

where $w$ satisfies

\begin{equation}\label{w stima lipschitz}
\vert{\:w}(x)\vert+\vert{\:x}\:\vert\:\vert{D}w(x)\vert\leq{C}\:\vert{\:x}\:\vert^
{\:2-n+\alpha},\quad{in}\quad{B}_{R}\setminus\{0\},
\end{equation}

\begin{equation}\label{w stima int}
\bigg( \int_{r<\vert{x}\vert<2r} \vert{D}^{2}w\vert^{p}\bigg) ^
{\frac{1}{p}}\leq{C}\:r^{\frac{n}{p}-n+\alpha},\quad\mbox{for}\:every\quad{r}
, \:0<r<R/2.
\end{equation}

Here $\alpha$ is such that $0<\alpha<1-\frac{n}{p}$, and C is a positive constant depending only on $\alpha,\:n,\:p,\:R$, $\lambda$, $E$, $\mathcal{E}$ and $k$.
\end{thm}

\begin{rem}\label{principal branch}
Since $K^{-1}(0)$ is a complex matrix, the expression 

\begin{equation}\label{u singular principal part}
\left(K^{-1}(0)x\cdot x\right)^{\frac{1}{2}} 
\end{equation}

appearing in the leading term in \eqref{u singular} is defined as the principal branch of \eqref{u singular principal part}, where a branch cut along the negative real axis of the complex plane has been defined for $z^{\frac{1}{2}}$, $z\in\mathbb{C}$. Expressions like \eqref{u singular principal part} will appear in the sequel of the paper and they will be understood in the same way.
\end{rem}

Next we consider two technical lemmas that are needed for the proof of Theorem \ref{theor singular sol}. The proofs of these results for the case where $L=-\mbox{div}(K\nabla\cdot)$, with $K$ a real matrix valued-function, are treated in detail in \cite {A1} and their extension to the more general case $L= - \mbox{div}(K\nabla\cdot)+q$, with $K,q$ a real matrix valued-function and a real function respectively, was extended in \cite{G}, therefore only the key points of their proof will be highlighted in the complex case below.


\begin{lem}\label{tech lemma 1}
Let $p>n$ and $u\in W^{2,p}_{loc}(B_R \setminus\{0\})$ be such that, for some positive $s$,

\begin{eqnarray}
& &|u(x)|\leq |x|^{2-s},\qquad\mbox{for any}\quad x\in B_R
\setminus\{0\},\\
& & \left(\int_{r<|x|<2r} |Lu|^p\right)^{\frac{1}{p}}\leq A
r^{\frac{n}{p}-s},\qquad\mbox{for\:any}\:r,\quad 0<r<\frac{R}{2}.
\end{eqnarray}

Then we have

\begin{eqnarray}
& & |Du(x)|\leq C |x|^{1-s},\quad\mbox{for any}\quad x\in
B_R\setminus\{0\},\label{estimate Du}\\
& &\left(\int_{r<|x|<2r} |D^2 u|^p\right)^{\frac{1}{p}}\leq C
r^{\frac{n}{p}-s}\quad\mbox{for any}\:r,\quad 0< r <
\frac{R}{4},\label{estimate D2u}
\end{eqnarray}
where $C$ is a positive constant depending only on $A$, $n$, $p$, $\lambda$, $E$, $\mathcal{E}$ and $k$.
\end{lem}

\textit{Proof of Lemma \ref{tech lemma 1}}. The proof of \eqref{estimate D2u} is based on the interior $L^{p}$ - Schauder estimate for uniformly elliptic systems

\begin{equation}\label{Schauder Lp}
\Big(\int_{r<|x|<2r} |D^2 u|^p\Big)^{\frac{1}{p}} \leq C \bigg\{\Big(\int_{\frac{r}{2}<|x|<4r} |L u|^p\Big)^{\frac{1}{p}}
 +  r^{-2} \Big( \int_{\frac{r}{2}<|x|<4r}|u|^p \Big)^{\frac{1}{p}}\bigg\},
\end{equation}



for every $r$, $0<r<\frac{R}{4}$, which, combined with interpolation inequality

\begin{eqnarray}\label{interpolation inequality}
r^{\frac{n}{p} -1}\sup_{r<|x|<2r} |Du(x)|  & \leq & C \Big\{\Big(\int_{\frac{r}{2}<|x|<4r} |D^2 u|^p\Big)^{\frac{1}{p}}\nonumber\\
& + & r^{-2} \Big( \int_{\frac{r}{2}<|x|<4r}|u|^p \Big)^{\frac{1}{p}}
\end{eqnarray}

leads to \eqref{estimate Du}. The positive constant $C$ appearing in \eqref{Schauder Lp} depends on $n$, $p$, $\lambda$, $E$, $\mathcal{E}$ and $k$ only, whereas the positive constant $C$ in \eqref{interpolation inequality} depends on $n$ and $p$ only. For \eqref{Schauder Lp} we refer to \cite[Lemma 6.2.6]{Mo}) and for a detailed proof of it, in the case of a single real equation in divergence form, we refer to \cite[Proof of Lemma 2.1]{A1}.  We refer to \cite[Theorem 5.12]{AF} for a detailed proof of \eqref{interpolation inequality} in the real case. For the complex case, \eqref{interpolation inequality} can be derived by denoting $u=u^1 + i u^2$ and combining 

\begin{eqnarray}\label{interpolation inequality complex}
& &r^{\frac{n}{p} -1}\sup_{r<|x|<2r} |Du^i(x)| \nonumber\\
& & \leq C \Big(||D^2 u^i||_{L^p (\frac{r}{2}<|x|<4r)} + r^{-2} ||u^i||_{L^p (\frac{r}{2}<|x|<4r)}\Big)\nonumber\\
& & \leq C \Big(||D^2 u||_{L^p (\frac{r}{2}<|x|<4r)} + r^{-2} ||u||_{L^p (\frac{r}{2}<|x|<4r)}\Big),
\end{eqnarray}

for $i=1,2$ together with

\begin{equation}
\sup_{r<|x|<2r} |Du(x)| \leq \sup_{r<|x|<2r} |Du^1(x)| + \sup_{r<|x|<2r} |Du^2(x)|. 
 \end{equation}
 
\hspace{11.5cm}$\blacksquare$


\begin{lem}\label{tech lemma 2}
Let $f\in L^{p}_{loc}(B_R \setminus\{0\})$ satisfy

\begin{equation}\label{estimate f}
\left(\int_{r<|x|<2r} |f|^p\right)^{\frac{1}{p}}\leq A
r^{\frac{n}{p}-s},\qquad\mbox{for any}\:r,\quad 0<r<\frac{R}{2},
\end{equation}

with $2<s<n<p$. Then there exists $u\in W^{2,p}_{loc}(B_R\setminus\{0\})$ satisfying

\begin{equation}\label{TLu}
Lu=f,\qquad\mbox{in}\quad B_R \setminus\{0\}
\end{equation}

and

\begin{equation}\label{estimate u}
|u(x)|\leq C |x|^{2-s},\quad\mbox{for any}\quad x\in
B_R\setminus\{0\},
\end{equation}

where $C$ is a positive constant depending only on $A$, $s$, $n$, $p$, $R$, $\lambda$, $E$, $\mathcal{E}$ and $k$.
\end{lem}

\textit{Proof of Lemma \ref{tech lemma 2}}. If $f\in L^{\infty}(B_R)$ then there exists a unique Green matrix $G(x,y)=\{G_{ij}(x,y)\}_{i,j=1}^{2}$ defined in $\{x,y\in B_R , x\neq y\}$ such that

\begin{equation}\label{green matrix}
LG(\cdot , y) = \delta(\cdot - y) I,\qquad\textnormal{for\: all}\: y\in B_R
\end{equation}

in the sense that for every $\phi = (\phi^1,\:\phi^2)\in C^{\infty}_c(B_R)$ we have

\begin{equation}
\int_{B_R} K^{\alpha\beta}_{ij} D_{\beta} G_{jk}(\cdot , y) D_{\alpha} \phi^i + q_{ij} G_{jk} (\cdot , y) \phi^i = \phi^k (y),\quad\textnormal{for}\: k=1,2.
\end{equation}

Moreover

\begin{equation}\label{fund sol complex}
|G (x,y)|\leq C |x-y|^{2-n},\qquad\mbox{for any}\:x\neq y,
\end{equation}

where $C$ is a positive constant depending on $n$, $\lambda$, $E$, $\mathcal{E}$ and $k$ and the vector valued - function $u=(u^1,\:u^2)$ defined by

\begin{equation}\label{u solution}
u^{k}(y)=\int_{B_R} G_{jk} (x,y) f^j(x) dx, \quad\textnormal{for}\: k=1,2,
\end{equation}

satisfies $Lu=f$ with

\begin{equation}\label{u bounds}
|u(x)|\leq \int_{B_R} |G(x,y)| |f(y)|\:dy\leq C(I_1 + I_2),
\end{equation}

where $f=(f^1,\:f^2)$ and

\begin{eqnarray}
I_1&=& \int_{|y|<\frac{|x|}{2}} |x-y|^{2-n} |f(y)|\:dy\label{I1}{\color{red}{,}}\\
I_2&=& \int_{\frac{|x|}{2}<|y|<R} |x-y|^{2-n} |f(y)|\:dy\label{I2}.
\end{eqnarray}

For the existence, uniqueness and asymptotic behaviour of the Green's matrix $G$ on $B_R$ as in \eqref{green matrix}-\eqref{fund sol complex}  we refer to \cite{DM}. We also refer to \cite{F}, \cite{F1} and the more recent result \cite{DaHMa} for further reading on the issue of the Green's matrix for elliptic systems of the second order. By an argument based on the monotone convergence theorem, one can show that $I_1$ and $I_2$ are both bounded from above by $C|x|^{2-s}$, where $C$ is a positive constant depending on $A,s,n,p,R,\lambda$, $E$, $\mathcal{E}$ and $k$.\\ 

If $f^{p}_{loc}(B_R\setminus\{0\})$,  we introduce a sequence $\{f_N\}_{N=1}^{\infty}$, with $f_N=(f^1_N,\:f^2_N)$, for $N\geq 1$, defined by

\begin{displaymath}\label{fN}
f^j_N = \left\{ \begin{array}{ll}  N
, &
\textrm{$\textnormal{when}\quad f^j>N$},\\
f^j &
\textrm{$\textnormal{when}\quad |f^j|\leq N$},\\
\!\!\!\!\!- N
, &
\textrm{$\textnormal{when}\quad f^j< -N$},
\end{array} \right.
\end{displaymath}

for $j=1,2$. $f_N\in L^{\infty} (B_R)$, for any $N\geq 1$ and $f_N\longrightarrow f$ pointwise on $B_R\setminus \{0\}$. For any $N\geq 1$, let $u_N\in W^{2,p}_{loc}(B_R\setminus \{0\})$ be the solution to 

\begin{equation}\label{uN 1}
Lu_N = f_N\qquad\textnormal{in}\quad B_R\setminus\{0\}
\end{equation}

such that 

\begin{equation}\label{uN 2}
|u_N(x)|\leq C_N |x|^{2-s},\qquad\textnormal{for\:any}\: x\in B_R\setminus\{0\}.
\end{equation}

$|f_N|\leq |f|$ on $B_R$, therefore $||f_N||_{L^p (\tilde\Omega)}\leq ||f||_{L^p (\tilde\Omega)}$, for any $\tilde\Omega$, $\tilde\Omega\subset\subset B_R\setminus\{0\}$, for any $N\geq 1$. By applying interior $L^p$ - Schauder estimates to $u_N$ and using the fact that $f\in L^{p}_{loc}(B_R\setminus\{0\})$ we obtain that 

\begin{equation}
||u_N||_{W^{2,p}(\tilde\Omega)}\leq C,\qquad\textnormal{for\:any}\:\tilde\Omega,\quad \tilde\Omega\subset\subset B_R\setminus\{0\},
\end{equation}

where $C$ is a positive constant that depends on $\tilde\Omega$. By applying a diagonal process we can find a subsequence $\{u_N\}_{N=1}^{\infty}$ weakly  converging in $W^{2,p}_{loc} (B_R\setminus\{0\})$ to some function $u\in W^{2,p}_{loc} (B_R\setminus\{0\})$. This limit satisfies both \eqref{TLu} and \eqref{estimate u}.\\ 

\hspace{11.5cm}$\blacksquare$\\


We proceed next with the proof of Theorem \ref{theor singular sol}.\\

\noindent\textit{Proof of Theorem \ref{theor singular sol}}  We start by considering 

\[H(x)=C\Big(K^{-1}(0) x\cdot x\Big)^{\frac{2-n}{2}},\]

solution to

\begin{equation}\label{H solution}
L_0 H = 0, \qquad\textnormal{in}\quad B_R\setminus\{0\},
\end{equation}

where $L_0:=- \mbox{div}\left(K(0)\nabla\cdot\right)$ on $B_R$. We want to find $w$ such  that

\begin{equation}\label{solution Lq}
L(H + w) = 0, \qquad\textnormal{in}\quad B_R\setminus\{0\},
\end{equation}

satisfying \eqref{w stima lipschitz}, \eqref{w stima int}, where $L$ is defined by \eqref{L complex}. We have

\begin{eqnarray}\label{Lq H}
-LH  & = & -L_0 H - L H \nonumber\\
& = & \Big(K_{ij} (x) - K_{ij}(0)\Big)\:\frac{\partial^{2}H}{\partial x_i \partial
x_j} - \frac{\partial a_{ij}}{\partial x_i}\:\frac{\partial
H}{\partial x_j} - qH.
\end{eqnarray}

Therefore for any $r$, $0<r<\frac{R}{2}$ we have

\begin{eqnarray}
\bigg(\int_{r<|x|<2r} |L H| ^p \bigg)^{\frac{1}{p}}  & \leq & \bigg(\int_{r<|x|<2r}
|K_{ij} (x) - K_{ij}(0)|^p\:\left|\frac{\partial^2 H}{\partial
x_i\partial x_j}\right|^p\bigg)^{\frac{1}{p}}\nonumber\\
& + & \bigg(\int_{r<|x|<2r} \left|\frac{\partial K_{ij}}{\partial
x_i}\right|^p\:\left|\frac{\partial H}{\partial
x_j}\right|^p\bigg)^{\frac{1}{p}}\nonumber\\
& + & \bigg(\int_{r<|x|<2r} |qH|^p \bigg)^{\frac{1}{p}}\nonumber\\
& \leq & \bigg(\int_{r<|x|<2r}
|x|^{\beta p}\:|x|^{-np}\bigg)^{\frac{1}{p}}\nonumber\\
& + & \bigg(\int_{r<|x|<2r} \left|\frac{\partial K_{ij}}{\partial
x_i}\right|^p\:|x|^{(1-n)p}\bigg)^{\frac{1}{p}}\nonumber\\
& + & \bigg(\lambda\int_{r<|x|<2r} |x|^{(2-n)p} \bigg)^{\frac{1}{p}} \nonumber\\
& \leq & C r^{\frac{n}{p} -n +\beta},
\end{eqnarray}

where $\beta = 1-\frac{n}{p}$ and $C$ is a positive constant depending on $\lambda$, $E$, $\mathcal{E}$, $R$ and $k$ only.  If we take $w\in W^{2,p}_{loc}(B_R\setminus\{0\})$ to be the solution to $Lw=f$ given by Lemma \ref{tech lemma 2}, with $f= -L H$ and $s=n-\beta$, then

\begin{equation}\label{estimate w}
|w(x)|\leq C |x|^{2-n+\beta}
\end{equation}

and, by Lemma \ref{tech lemma 1}, properties \eqref{w stima lipschitz}, \eqref{w stima int} are satisfied. $\hspace{3cm}\blacksquare$\\




\section{Proof of the main result.}\label{proofs main result}

Since the boundary $\partial\Omega$ is Lipschitz, the normal unit vector field might not be defined on $\partial\Omega$. We shall therefore introduce a unitary vector field $\widetilde\nu$ locally defined near $\partial\Omega$ such that: (i) $\widetilde\nu$ is $C^{\infty}$ smooth, (ii) $\widetilde\nu$ is non-tangential to $\partial\Omega$ and it points to the exterior of $\Omega$
(see \cite[Lemmas 3.1-3.3]{AG} for a precised construction of $\widetilde\nu$). Here we simply recall that any point $z_{\tau}=x^{0}+\tau\widetilde\nu$, where $x^{0}\in\partial\Omega$, satisfies

\begin{equation}\label{def tau0}
C\:\tau\leq{d}(z_{\tau},\:\partial\Omega)\leq\tau,
\quad\textnormal{for}\:\textnormal{any}\quad\tau,\quad
0\leq\tau\leq\tau_{0},
\end{equation}

where $\tau_{0}$ and $C$ depend on $L$, $r_0$ only.


\begin{rem}\label{remark C}
Several constants depending on the \textit{a-priori data} introduced in Definition \ref{a priori data} will appear in the proof of the main result below. In order to simplify our notation, we shall denote by $C$ any of these constants, avoiding in most cases to point out their specific dependence on the \textit{a-priori data} which may vary from case to case. 
\end{rem}


\textit{Proof of Theorem \ref{stabilita' al bordo}.}  We start by recalling that by \eqref{Alessandrini identity} we have

\begin{eqnarray*}
\langle (\Lambda_{\mu_{a_1}} - \Lambda_{\mu_{a_2}} ) u,\overline{v}\rangle &=& \int_{\Omega} \left(K_{\mu_{a_1}}(x) - K_{\mu_{a_2}}(x)\right)\nabla u(x)\cdot\nabla v(x)\:dx\nonumber\\
&+&\int_{\Omega}\left(\mu_{a_1}(x)-\mu_{a_2}(x)\right)u(x)v(x)\:dx,
\end{eqnarray*}

for any $u,v \in H^1 (\Omega)$ that solve

\begin{eqnarray}
& & \mbox{div}\big(K_{\mu_{a_1}} \nabla u\big) + (\mu_{a_1} -ik)u = 0,\qquad\textnormal{in}\quad\Omega ,\label{eq1}\\
& & \mbox{div}\big(K_{\mu_{a_2}} \nabla v\big) + (\mu_{a_2} -ik)v = 0,\qquad\textnormal{in}\quad\Omega .\label{eq2}
\end{eqnarray}

We set $x^{0}\in\partial\Omega$ such that

\begin{equation*}
(\mu_{a_1}-\mu_{a_2})(x^{0})\:=\:\parallel{\mu_{a_1}}-\mu_{a_2}\parallel_{L^{\infty}(\partial\Omega)}
\end{equation*}

and $z_{\tau}=x^{0}+\tau\widetilde\nu$, with $0<\tau\leq\tau_{0}$, where $\tau_{0}$ is the number fixed in \eqref{def tau0}. Let $u,\:v\in W^{2,p}(\Omega)$ be the singular solutions of Theorem \ref{theor singular sol} to \eqref{eq1}, \eqref{eq2} respectively, having a singularity at $z_{\tau}$ 

\begin{eqnarray}\label{singular ui}
u(x) & =& \Big(K^{-1}_{\mu_{a_1}}(z_{\tau})(x-z_{\tau})\cdot(x-z_{\tau})\Big)^{\frac{2-n}{2}}+ O\left(\left|x-z_{\tau}\right|^{2-n+\alpha}\right),\nonumber\\
v(x) &=& \Big(K^{-1}_{\mu_{a_2}}(z_{\tau})(x-z_{\tau})\cdot(x-z_{\tau})\Big)^{\frac{2-n}{2}}+ O\left(\left|x-z_{\tau}\right|^{2-n+\alpha}\right).
\end{eqnarray}

By setting $\rho=2\tau_{0}$ we have that $B_{\rho}(z_{\tau})\cap\Omega\neq\emptyset$ and from \eqref{Alessandrini identity} we obtain

\begin{eqnarray}\label{integral inequality complex 1}
& &\parallel\Lambda_{\mu_{a_1}}-\Lambda_{{\mu_{a_2}}}\parallel_{\mathcal{L}(H^{\frac{1}{2}}(\partial\Omega),H^{-\frac{1}{2}}(\partial\Omega))}\||\overline{u}||_{H^{\frac{1}{2}}(\partial\Omega)}\:||v||_{H^{\frac{1}{2}}(\partial\Omega)}\nonumber\\
&  &\geq  \left|\int_{\Omega\cap{B}_{\rho}(z_{\tau})}\Big(K_{\mu_{a_1}}(x)-K_{\mu_{a_2}}(x)\Big)\:\nabla u(x)\cdot\nabla v(x)\:dx\right| \nonumber\\
& & - \int_{\Omega\setminus{B}_{\rho}(z_{\tau})}\Big|K_{\mu_{a_1}}(x)-K_{\mu_{a_2}}(x)\Big|\:|\nabla u(x)|\:|\nabla v(x)|\:dx\nonumber\\
& & - \int_{\Omega\cap{B}_{\rho}(z_{\tau})}\Big|(\mu_{a_1}-\mu_{a_2})(x)\Big|\left|u(x)\right|\!\left|v(x)\right|\:dx\nonumber\\
& & - \int_{\Omega\setminus{B}_{\rho}(z_{\tau})}\Big|(\mu_{a_1}-\mu_{a_2})(x)\Big|\left|u(x)\right|\!\left|v(x)\right|\:dx.
\end{eqnarray}

By \eqref{singular ui} and Theorem \ref{theor singular sol} we have 

\begin{eqnarray}\label{grad singular ui}
\nabla u(x) &=& (2-n)\big(K^{-1}_{\mu_{a_1}}(z_{\tau})(x-z_{\tau})\cdot (x-z_{\tau})\big)^{-\frac{n}{2}}\: K^{-1}_{\mu_{a_1}}(z_{\tau}) (x-z_{\tau})\nonumber\\
& +& O(|x-z_{\tau}|^{1-n+\alpha}),\nonumber\\
\nabla v(x)&=&(2-n)\big(K^{-1}_{\mu_{a_2}}(z_{\tau})(x-z_{\tau})\cdot (x-z_{\tau})\big)^{-\frac{n}{2}}\: K^{-1}_{\mu_{a_2}}(z_{\tau}) (x-z_{\tau}) \nonumber\\
&+& O(|x-z_{\tau}|^{1-n+\alpha}).
\end{eqnarray}

Recalling that for $i=1,2$ the real and imaginary parts of $K^{-1}_{\mu_{a_i}}$ satisfy \eqref{apriori ass ineq 1} and \eqref{apriori ass ineq 2} respectively, we have

\begin{equation}\label{K inverse module estimate}
C^{-1}|\xi|^2\leq \big| K^{-1}_{\mu_{a_i}}(x) \xi\cdot\xi\big| \leq C |\xi|^2,\quad\textnormal{for\:a.e.}\: x\in\Omega,\quad\textnormal{for\:every}\:\xi\in\mathbb{R}^n
\end{equation}

and by combining \eqref{integral inequality complex 1} together with \eqref{singular ui}, \eqref{grad singular ui} and \eqref{K inverse module estimate} we obtain

\begin{eqnarray}\label{integral inequality complex 2}
& & \left|\int_{\Omega\cap{B}_{\rho}(z_{\tau})}\Big(K_{\mu_{a_1}}(x)-K_{\mu_{a_2}}(x)\Big)\:\nabla u(x)\cdot\nabla v(x)\:dx\right|\nonumber\\
& &\leq\:C\bigg\{\int_{\Omega\cap{B}_{\rho}(z_{\tau})}\:\vert x-\:z_{\tau}\vert^{4-2n}\:dx\nonumber\\
& &+ \int_{\Omega\setminus B_{\rho}(z_{\tau})}\:\vert x-\:z_{\tau}\vert^{4-2n}\:dx\nonumber\\
& &\:+\:\int_{\Omega\setminus{B}_{\rho}(z_{\tau})}\:\vert x-\:z_{\tau}\vert^{2-2n}\:dx\nonumber\\
& & + \parallel\Lambda_{\mu_{a_1}}-\Lambda_{{\mu_{a_2}}}\parallel_{\mathcal{L}(H^{\frac{1}{2}}(\partial\Omega),H^{-\frac{1}{2}}(\partial\Omega))}\||u||_{H^{\frac{1}{2}}(\partial\Omega)}\:||v||_{H^{\frac{1}{2}}(\partial\Omega)}\bigg\}.
\end{eqnarray}

The left hand side of \eqref{integral inequality complex 2} can be estimated from below by recalling that $K_{\mu_{a_i}}(\cdot)$ is H\"older continuous on $\overline\Omega$ with exponent $\beta=1-\frac{n}{p}$, for $i=1,2$ and by recalling again \eqref{singular ui}, which leads to

\begin{eqnarray}\label{integral inequality complex 3}
& & \left|\int_{\Omega\cap{B}_{\rho}(z_{\tau})}\Big(K_{\mu_{a_1}}(x)-K_{\mu_{a_2}}(x)\Big)\:\nabla u(x)\cdot\nabla v(x)\:dx\right|\nonumber\\
& & \geq  \left|\int_{\Omega\cap{B}_{\rho}(z_{\tau})}\Big(K_{\mu_{a_1}}(x^0)-K_{\mu_{a_2}}(x^0)\Big)\:\nabla u(x)\cdot\nabla v(x)\:dx\right|\nonumber\\
& & -  C \int_{\Omega\cap{B}_{\rho}(z_{\tau})} |x-x^0|^{\beta}\:|\:\nabla u(x)|\: |\nabla v(x)|\:dx\nonumber\\
& & \geq  \left|\int_{\Omega\cap{B}_{\rho}(z_{\tau})}\Big(K_{\mu_{a_1}}(x^0)-K_{\mu_{a_2}}(x^0)\Big)\:\nabla u(x)\cdot\nabla v(x)\:dx\right|\nonumber\\
& & -  C \int_{\Omega\cap{B}_{\rho}(z_{\tau})} |x-x^0|^{\beta}\:|x-z_{\tau}|^{2-2n}\:dx
\end{eqnarray}

and by combining \eqref{integral inequality complex 3} together with \eqref{integral inequality complex 2} we obtain

\begin{eqnarray}\label{integral inequality complex 4}
& & \left|\int_{\Omega\cap{B}_{\rho}(z_{\tau})}\Big(K_{\mu_{a_1}}(x^0)-K_{\mu_{a_2}}(x^0)\Big)\:\nabla u(x)\cdot\nabla v(x)\:dx\right|\nonumber\\
& &\leq \:C\bigg\{ \int_{\Omega\cap{B}_{\rho}(z_{\tau})}\:\vert x-\:z_{\tau}\vert^{2-2n}\: |x-x^0|^{\beta}\:dx\nonumber\\
& &+ \int_{\Omega\cap{B}_{\rho}(z_{\tau})}\:\vert x-\:z_{\tau}\vert^{4-2n}\:dx\nonumber\\
& &+ \int_{\Omega\setminus B_{\rho}(z_{\tau})}\:\vert x-\:z_{\tau}\vert^{4-2n}\:dx\nonumber\\
& &\:+\:\int_{\Omega\setminus{B}_{\rho}(z_{\tau})}\:\vert x-\:z_{\tau}\vert^{2-2n}\:dx\nonumber\\
& & + \parallel\Lambda_{\mu_{a_1}}-\Lambda_{{\mu_{a_2}}}\parallel_{\mathcal{L}(H^{\frac{1}{2}}(\partial\Omega),H^{-\frac{1}{2}}(\partial\Omega))}\||u||_{H^{\frac{1}{2}}(\partial\Omega)}\:||v||_{H^{\frac{1}{2}}(\partial\Omega)}\bigg\}.
\end{eqnarray}

Recalling \eqref{K inverse module estimate} and combining it together with \eqref{grad singular ui}, we can estimate the left hand side of \eqref{integral inequality complex 4} from below as

\begin{eqnarray}\label{integral inequality complex 5}
& & \left|\int_{\Omega\cap{B}_{\rho}(z_{\tau})}\Big(K_{\mu_{a_1}}(x^0)-K_{\mu_{a_2}}(x^0)\Big)\:\nabla u(x)\cdot\nabla v(x)\:dx\right|\nonumber\\
& &\geq  (2-n)^2  \nonumber\\
& & \times \left|\int_{\Omega\cap{B}_{\rho}(z_{\tau})}\frac{K^{-1}_{\mu_{a_2}}(z_{\tau})\Big(K_{\mu_{a_1}}(x^0)-K_{\mu_{a_2}}(x^0)\Big)K^{-1}_{\mu_{a_1}}(z_{\tau})\:(x-z_{\tau})\cdot (x-z_{\tau})}{\big(K^{-1}_{\mu_{a_1}}(z_{\tau})(x-z_{\tau})\cdot (x-z_{\tau})\big)^{\frac{n}{2}} \big(K^{-1}_{\mu_{a_2}}(z_{\tau})(x-z_{\tau})\cdot (x-z_{\tau})\big)^{\frac{n}{2}}}\:dx\right|\nonumber\\
& & -C\bigg\{ \int_{\Omega\cap{B}_{\rho}(z_{\tau})} |x-z_{\tau}|^{2-2n+\alpha} +  \int_{\Omega\cap{B}_{\rho}(z_{\tau})} |x-z_{\tau}|^{2-2n+2\alpha}\bigg\}.
\end{eqnarray}

\eqref{integral inequality complex 5} together with \eqref{integral inequality complex 4} leads to

\begin{eqnarray}\label{integral inequality complex 6}
& & \left|\int_{\Omega\cap{B}_{\rho}(z_{\tau})}\frac{K^{-1}_{\mu_{a_2}}(z_{\tau})\Big(K_{\mu_{a_1}}(x^0)-K_{\mu_{a_2}}(x^0)\Big)K^{-1}_{\mu_{a_1}}(z_{\tau})\:(x-z_{\tau})\cdot (x-z_{\tau})}{\big(K^{-1}_{\mu_{a_1}}(z_{\tau})(x-z_{\tau})\cdot (x-z_{\tau})\big)^{\frac{n}{2}} \big(K^{-1}_{\mu_{a_2}}(z_{\tau})(x-z_{\tau})\cdot (x-z_{\tau})\big)^{\frac{n}{2}}}\:dx\right|\nonumber\\
& &\leq \:C\bigg\{ \int_{\Omega\cap{B}_{\rho}(z_{\tau})}\:\vert x-\:z_{\tau}\vert^{2-2n+\alpha}\:dx + \int_{\Omega\cap{B}_{\rho}(z_{\tau})}\:\vert x-\:z_{\tau}\vert^{2-2n}\: |x-x^0|^{\beta}\:dx\nonumber\\
& &+ \int_{\Omega\cap{B}_{\rho}(z_{\tau})}\:\vert x-\:z_{\tau}\vert^{4-2n}\:dx + \int_{\Omega\setminus B_{\rho}(z_{\tau})}\:\vert x-\:z_{\tau}\vert^{4-2n}\:dx\nonumber\\
& & +\:\int_{\Omega\setminus{B}_{\rho}(z_{\tau})}\:\vert x-\:z_{\tau}\vert^{2-2n}\:dx \nonumber\\
& & + \parallel\Lambda_{\mu_{a_1}}-\Lambda_{{\mu_{a_2}}}\parallel_{\mathcal{L}(H^{\frac{1}{2}}(\partial\Omega),H^{-\frac{1}{2}}(\partial\Omega))}\||u||_{H^{\frac{1}{2}}(\partial\Omega)}\:||v||_{H^{\frac{1}{2}}(\partial\Omega)}\bigg\}.
\end{eqnarray}

$K^{-1}_{\mu_{a_i}}$ is H\"older continuous on $\overline\Omega$, with $\beta=1-\frac{n}{p}$, for $i=1,2$ and, recalling that $C\tau\leq |x-z_{\tau}|$, we have

\begin{eqnarray}\label{estimate numerator}
& & K^{-1}_{\mu_{a_2}}(z_{\tau}))\Big(K_{\mu_{a_1}}(x^0)-K_{\mu_{a_2}}(x^0)\Big)K^{-1}_{\mu_{a_1}}(z_{\tau})\:(x-z_{\tau})\cdot (x-z_{\tau})\nonumber\\
& & =\Big(K^{-1}_{\mu_{a_2}}(x^0) + O(\tau^{\beta})\Big)\Big(K_{\mu_{a_1}}(x^0)-K_{\mu_{a_2}}(x^0)\Big)\Big( K^{-1}_{\mu_{a_1}}(x^0 )+ O(\tau^{\beta})\Big)\:(x-z_{\tau})\cdot (x-z_{\tau})\nonumber\\
& & = \Big(K^{-1}_{\mu_{a_2}}(x^0) - K^{-1}_{\mu_{a_1}}(x^0)\Big)\:(x-z_{\tau})\cdot (x-z_{\tau})+ O(|x-z_{\tau}|^{2+\beta})\nonumber\\
& & =n (\mu_{a_2 }- \mu_{a_1})(x^0)\:|x-z_{\tau}|^2+ O(|x-z_{\tau}|^{2+\beta}).
\end{eqnarray}

Hence \eqref{integral inequality complex 6}, combined with \eqref{estimate numerator} and again with \eqref{K inverse module estimate}, leads to

\begin{eqnarray}\label{integral inequality complex 7}
& &(\mu_{a_1 }- \mu_{a_2})(x^0)\nonumber\\
& & \times\left|\int_{\Omega\cap{B}_{\rho}(z_{\tau})}\frac{|x-z_{\tau}|^2}{\big(K^{-1}_{\mu_{a_1}}(z_{\tau})(x-z_{\tau})\cdot (x-z_{\tau})\big)^{\frac{n}{2}} \big(K^{-1}_{\mu_{a_2}}(z_{\tau})(x-z_{\tau})\cdot (x-z_{\tau})\big)^{\frac{n}{2}}}\:dx\right|\nonumber\\
& &\leq \:C\bigg\{ \int_{\Omega\cap{B}_{\rho}(z_{\tau})}\:\vert x-\:z_{\tau}\vert^{2-2n+\beta}\:dx +\int_{\Omega\cap{B}_{\rho}(z_{\tau})}\:\vert x-\:z_{\tau}\vert^{2-2n+\alpha}\:dx\nonumber\\
& & + \int_{\Omega\cap{B}_{\rho}(z_{\tau})}\:\vert x-\:z_{\tau}\vert^{2-2n}\: |x-x^0|^{\beta}\:dx + \int_{\Omega\cap{B}_{\rho}(z_{\tau})}\:\vert x-\:z_{\tau}\vert^{4-2n}\:dx\nonumber\\
& &+ \int_{\Omega\setminus B_{\rho}(z_{\tau})}\:\vert x-\:z_{\tau}\vert^{4-2n}\:dx +\:\int_{\Omega\setminus{B}_{\rho}(z_{\tau})}\:\vert x-\:z_{\tau}\vert^{2-2n}\:dx\nonumber\\
& & + \parallel\Lambda_{\mu_{a_1}}-\Lambda_{{\mu_{a_2}}}\parallel_{\mathcal{L}(H^{\frac{1}{2}}(\partial\Omega),H^{-\frac{1}{2}}(\partial\Omega))}\:\||u||_{H^{\frac{1}{2}}(\partial\Omega)}\:||v||_{H^{\frac{1}{2}}(\partial\Omega)}\bigg\}.
\end{eqnarray}

The integrand appearing on the left hand side of \eqref{integral inequality complex 7} can be expressed as

\begin{equation}\label{integrand 7}
\frac{|x-z_{\tau}|^2\: F(x)}{\big|K^{-1}_{\mu_{a_1}}(z_{\tau})(x-z_{\tau})\cdot (x-z_{\tau})\big|^{n} \big|K^{-1}_{\mu_{a_2}}(z_{\tau})(x-z_{\tau})\cdot (x-z_{\tau})\big|^{n}},
\end{equation}

where the complex-valued function $F$ is defined by

\begin{equation}\label{F}
F(x):=\bigg\{\Big(\overline{K}^{-1}_{\mu_{a_1}} (z_{\tau})(x-z_{\tau})\cdot(x-z_{\tau})\Big)\Big(\overline{K}^{-1}_{\mu_{a_2}} (z_{\tau})(x-z_{\tau})\cdot(x-z_{\tau})\Big)\bigg\}^{\frac{n}{2}}.
\end{equation}

The choices of $k$ in either  \eqref{k range 1} or \eqref{k range 2} imply 

\begin{equation}\label{conditions on F}
|\Im F(x)| \leq |\Re F(x)|\quad\textnormal{and}\quad \Re F(x)>0{\color{red}{,}}
\end{equation}

where $\Re z$ and $\Im z$ denote the real and imaginary parts of a complex number $z$ respectively. By combining \eqref{conditions on F} together with \eqref{K inverse module estimate},  the left hand side of inequality \eqref{integral inequality complex 7} can be estimated from below as

\begin{eqnarray}\label{integral inequality complex 8}
& &(\mu_{a_1 }- \mu_{a_2})(x^0)\nonumber\\
& & \times\left|\int_{\Omega\cap{B}_{\rho}(z_{\tau})}\frac{|x-z_{\tau}|^2\: F(x)}{\big|K^{-1}_{\mu_{a_1}}(z_{\tau})(x-z_{\tau})\cdot (x-z_{\tau})\big|^{n} \big|K^{-1}_{\mu_{a_2}}(z_{\tau})(x-z_{\tau})\cdot (x-z_{\tau})\big|^{n}}\:dx\right|\nonumber\\
& &\geq (\mu_{a_1 }- \mu_{a_2})(x^0) \nonumber\\
& & \times \Re \bigg[\int_{\Omega\cap{B}_{\rho}(z_{\tau})}\frac{|x-z_{\tau}|^2\: F(x)}{\big|K^{-1}_{\mu_{a_1}}(z_{\tau})(x-z_{\tau})\cdot (x-z_{\tau})\big|^{n} \big|K^{-1}_{\mu_{a_2}}(z_{\tau})(x-z_{\tau})\cdot (x-z_{\tau})\big|^{n}}\:dx\bigg]\nonumber\\
& & \geq\frac{1}{\sqrt{2}} (\mu_{a_1 }- \mu_{a_2})(x^0)\nonumber\\
& &\times \int_{\Omega\cap{B}_{\rho}(z_{\tau})}\frac{|x-z_{\tau}|^2\: |F(x)|}{\big|K^{-1}_{\mu_{a_1}}(z_{\tau})(x-z_{\tau})\cdot (x-z_{\tau})\big|^{n} \big|K^{-1}_{\mu_{a_2}}(z_{\tau})(x-z_{\tau})\cdot (x-z_{\tau})\big|^{n}}\:dx\nonumber\\
& & \geq\frac{1}{\sqrt{2}} (\mu_{a_1 }- \mu_{a_2})(x^0)\nonumber\\
& & \times\!\!\int_{\Omega\cap{B}_{\rho}(z_{\tau})}\hspace{-0.85cm}\frac{|x-z_{\tau}|^2\: \Big|\overline{K}^{-1}_{\mu_{a_1}} (z_{\tau})(x-z_{\tau})\cdot(x-z_{\tau})\Big|^{\frac{n}{2}}\:\Big|\overline{K}^{-1}_{\mu_{a_2}} (z_{\tau})(x-z_{\tau})\cdot(x-z_{\tau})\Big|^{\frac{n}{2}}}{\big|K^{-1}_{\mu_{a_1}}(z_{\tau})(x-z_{\tau})\cdot (x-z_{\tau})\big|^{n} \big|K^{-1}_{\mu_{a_2}}(z_{\tau})(x-z_{\tau})\cdot (x-z_{\tau})\big|^{n}}\:dx.\nonumber\\
\end{eqnarray}

Combing \eqref{integral inequality complex 8} together with \eqref{K inverse module estimate}, we obtain

\begin{eqnarray}\label{integral inequality complex 8 bis}
& &(\mu_{a_1 }- \mu_{a_2})(x^0)\nonumber\\
& & \times\left|\int_{\Omega\cap{B}_{\rho}(z_{\tau})}\frac{|x-z_{\tau}|^2\: F(x)}{\big|K^{-1}_{\mu_{a_1}}(z_{\tau})(x-z_{\tau})\cdot (x-z_{\tau})\big|^{n} \big|K^{-1}_{\mu_{a_2}}(z_{\tau})(x-z_{\tau})\cdot (x-z_{\tau})\big|^{n}}\:dx\right|\nonumber\\
& & \geq\frac{1}{\sqrt{2}} (\mu_{a_1 }- \mu_{a_2})(x^0)\:C\int_{\Omega\cap{B}_{\rho}(z_{\tau})} |x-z_{\tau}|^{2-2n}\:dx.
\end{eqnarray}

\eqref{integral inequality complex 8 bis} combined with \eqref{integral inequality complex 7} and \eqref {integrand 7} then leads to

\begin{eqnarray}\label{integral inequality complex 9}
& &||\mu_{a_1 }- \mu_{a_2}||_{L^{\infty}(\partial\Omega)}\int_{\Omega\cap{B}_{\rho}(z_{\tau})} |x-z_{\tau}|^{2-2n}\:dx.\nonumber\\
& &\leq \:C\bigg\{ \int_{\Omega\cap{B}_{\rho}(z_{\tau})}\:\vert x-\:z_{\tau}\vert^{2-2n+\beta}\:dx +\int_{\Omega\cap{B}_{\rho}(z_{\tau})}\:\vert x-\:z_{\tau}\vert^{2-2n+\alpha}\:dx\nonumber\\
& & + \int_{\Omega\cap{B}_{\rho}(z_{\tau})}\:\vert x-\:z_{\tau}\vert^{2-2n}\: |x-x^0|^{\beta}\:dx + \int_{\Omega\cap{B}_{\rho}(z_{\tau})}\:\vert x-\:z_{\tau}\vert^{4-2n}\:dx\nonumber\\
& &+ \int_{\Omega\setminus B_{\rho}(z_{\tau})}\:\vert x-\:z_{\tau}\vert^{4-2n}\:dx +\:\int_{\Omega\setminus{B}_{\rho}(z_{\tau})}\:\vert x-\:z_{\tau}\vert^{2-2n}\:dx\nonumber\\
& & + \parallel\Lambda_{\mu_{a_1}}-\Lambda_{{\mu_{a_2}}}\parallel_{\mathcal{L}(H^{\frac{1}{2}}(\partial\Omega),H^{-\frac{1}{2}}(\partial\Omega))}\||u||_{H^{\frac{1}{2}}(\partial\Omega)}\:||v||_{H^{\frac{1}{2}}(\partial\Omega)}\bigg\}.
\end{eqnarray}

By recalling \eqref{def tau0}, the first integral appearing on the right hand side of \eqref{integral inequality complex 9} can be estimated from above by observing that $\Omega\cap{B}_{\rho}(z_{\tau})\subset\left\{x\:|\: C\tau\leq|x-z_{\tau}|\leq 2\tau_0\right\}$, therefore

\begin{eqnarray}\label{I1}
\int_{\Omega\cap{B}_{\rho}(z_{\tau})}\:\vert x-\:z_{\tau}\vert^{2-2n+\beta}\:dx &\leq &  \int_{\left\{C\tau\leq |x-z_{\tau}|\leq 2\tau_0\right\}}\:\vert x-\:z_{\tau}\vert^{2-2n+\beta}\:dx\nonumber\\
&=& \int_{C\tau}^{2\tau_0} s^{2-2n +\beta +n -1}\: ds  \int_{\{ |\xi | = 1\}} dS_{\xi}\nonumber\\
& \leq & C\left((C\tau)^{2-n+\beta} - (2\tau_0)^{2-n+\beta}\right)\nonumber\\
& \leq & C\tau^{2-n+\beta},
\end{eqnarray}

(see also \cite{A1}, \cite{AG}), where $dS_{\xi}$ denotes the surface measure on the unit sphere. Similarly to \eqref{I1}, the second, third and forth integrals on the right hand side of inequality \eqref{integral inequality complex 9} are estimated from above as

\begin{eqnarray}\label{I2 I3 I4}
\int_{\Omega\cap{B}_{\rho}(z_{\tau})}\:\vert x-\:z_{\tau}\vert^{2-2n+\alpha}\:dx &\leq & C\tau^{2-n+\alpha},\nonumber\\
\int_{\Omega\cap{B}_{\rho}(z_{\tau})}\:\vert x-\:z_{\tau}\vert^{2-2n}\:|x-x^0|^{\beta}\:dx &\leq & C\tau^{2-n+\beta},\nonumber\\
\int_{\Omega\cap{B}_{\rho}(z_{\tau})}\:\vert x-\:z_{\tau}\vert^{4-2n}\:dx &\leq & C\tau^{4-n}.
\end{eqnarray}

By observing that $\left(\Omega\setminus B_{\rho}(z_{\tau})\right)\subset \left\{x\:|\: 2\tau_0\leq|x-z_{\tau}|\leq R\right\}$, where $R$ depends on $\textnormal{diam}(\Omega)$, the last two integrals appearing on the right hand side of \eqref{integral inequality complex 9} can be estimated from above as

\begin{eqnarray}\label{I5 I6}
\int_{\Omega\setminus B_{\rho}(z_{\tau})}\:\vert x-\:z_{\tau}\vert^{4-2n}\:dx &\leq & \int_{\left\{2\tau_0\leq |x-z_{\tau}|\leq R \right\}}\:\vert x-\:z_{\tau}\vert^{4-2n}\:dx \leq C, \nonumber\\
\int_{\Omega\setminus B_{\rho}(z_{\tau})}\:\vert x-\:z_{\tau}\vert^{2-2n}\:dx &\leq & C.
\end{eqnarray}

The integral appearing on the left hand side of \eqref{integral inequality complex 9} can be estimated from below as

\begin{equation}\label{I0 from below}
\int_{\Omega\cap{B}_{\rho}(z_{\tau})} |x-z_{\tau}|^{2-2n}\:dx\geq C\tau^{2-n}
\end{equation}

and we refer to \cite[p.66]{Sa} for a detailed calculation of estimate \eqref{I0 from below}. By combining \eqref {integral inequality complex 9} together with \eqref{I1} - \eqref{I0 from below} and the $H^{\frac{1}{2}}(\partial\Omega)$ norms of {\color{red}{$u, v$}} (see \cite{A1}, \cite{AG}), we obtain

\begin{eqnarray}\label{inequality tau 1}
& & \parallel\mu_{a_1}-\mu_{a_2}\parallel_{L^{\infty}(\partial\Omega)}\tau^{2-n}  \leq C  \Big\{ \tau^{2-n+\beta}+\tau^{2-n+\alpha}+\tau^{4-n}+C\nonumber\\
& & +   \tau^{2-n}\parallel\Lambda_{\mu_{a_1}}-\Lambda_{{\mu_{a_2}}}\parallel_{\mathcal{L}(H^{\frac{1}{2}}(\partial\Omega),H^{-\frac{1}{2}}(\partial\Omega))}\Big\}.
\end{eqnarray}

By multiplying \eqref{inequality tau 1} by $\tau^{n-2}$ we obtain

\begin{equation}\label{inequality tau 2}
\parallel\mu_{a_1}-\mu_{a_2}\parallel_{L^{\infty}(\partial\Omega)}\leq C
\left\{\omega(\tau)+\parallel\Lambda_{\mu_{a_1}}-\Lambda_{{\mu_{a_2}}}\parallel_{\mathcal{L}(H^{\frac{1}{2}}(\partial\Omega),H^{-\frac{1}{2}}(\partial\Omega))}\right\},
\end{equation}

where $\omega(\tau)\rightarrow 0$ as $\tau\rightarrow{0}$, which concludes the proof. $\quad\blacksquare$\\

\begin{rem}
When $n=3$ the ranges for $k$, \eqref{k range 1} and \eqref{k range 2}, simplify to

\begin{equation}\label{k range 1 3D}
0<k \leq k_0 := \sqrt{3 \lambda^2 (1+\mathcal{E})^2 + \lambda^{-2} (1+\mathcal{E}^{-1})^{2}} - \sqrt{3}\lambda (1+\mathcal{E}),
\end{equation}

and

\begin{equation}\label{k range 2 3D}
k \geq \tilde{k}_0 :=(2+\sqrt{3})\: \lambda (1+\mathcal{E}).
\end{equation}

\end{rem}

\section*{\normalsize{Acknowledgements}}
R. Gaburro would like to express her gratitude to Giovanni Alessandrini for the fruitful conversations regarding the construction of the singular solutions in this manuscript and to Simon Arridge who kindly offered insights on the physical problem of diffuse Optical Tomography some time ago. R. Gaburro would like to acknowledge the support of the Higher Education Authority (HEA) Government of Ireland International Academic Mobility Program 2019. W. Lionheart would like to acknowledge the support of EPSRC grants EP/F033974/1 and EP/L019108/1 and the Royal Society for a Wolfson Research Merit Award.



\begin{thebibliography}{AB}

\bibitem{AF} R. Adams and J. Fournier, Sobolev Spaces, Academic Press, Amsterdam (2003).



\bibitem{A1} G. Alessandrini, Singular solutions of elliptic equations and the determination of conductivity by boundary measurements, J.
Differential Equations \textbf{84}, (2) (1990), 252-272.



\bibitem{A-dH-F-G-S} G. Alessandrini, F. Faucher, M. V. de Hoop, R. Gaburro and E. Sincich, Inverse problem for the Helmholtz equation with Cauchy data: reconstruction with conditional well-posedness driven iterative regularization, ESAIM: Mathematical Modeliing and Numerical Analysis \textbf{53} (\textbf{3}) (2019), 1005-1030.


\bibitem{A-dH-G} G. Alessandrini, M.V. de Hoop and R. Gaburro, Uniqueness for the electrostatic inverse boundary value problem with piecewise constant anisotropic conductivities, Inverse Problems \textbf{33} (\textbf{12}) (2017), 125013.

\bibitem{A-dH-G-S} G. Alessandrini, M. V. de Hoop, R. Gaburro and E. Sincich, Lipschitz stability for the electrostatic inverse boundary value problem with piecewise linear conductivities, Journal de Math\`ematiques Pures et Appliqu\`ees \textbf{107} (\textbf{5}) (2017), 638 - 664.

\bibitem{A-dH-G-S1} G. Alessandrini, M. V. de Hoop, R. Gaburro and E. Sincich, Lipschitz stability for a piecewise linear Schr\"odinger potential from local Cauchy data, Asymptotic Analysis \textbf{108} (\textbf{3}) (2018), 115-149.

\bibitem{A-dH-G-S2} G. Alessandrini, M. V. de Hoop, R. Gaburro and E. Sincich, EIT in a layered anisotropic medium, Inverse Problems and Imaging \textbf{12} (\textbf{3}) (2018), 667 - 676.

\bibitem{AG} G. Alessandrini and R. Gaburro, Determining conductivity with special anisotropy by boundary measurements, SIAM J. MATH. ANAL. \textbf{33}  (2001), 153-171.

\bibitem{AG1} G. Alessandrini and R. Gaburro, The local Calder\'on problem and the determination at the boundary of the conductivity, Comm. Partial
Differential Equations \textbf{34}  (2009), 918-936.


\bibitem{A-V} G. Alessandrini and S. Vessella, Lipschitz stability for the inverse conductivity problem, Advances in Applied Mathematics \textbf{35} (2005), 207-241.

\bibitem{Ar} S. R. Arridge, Optical tomography in medical imaging, Inverse Problems \textbf{15} (\textbf{2}) (1999), R41.

\bibitem{HebdenArridge} S. R. Arridge and J. C. Hebden, Optical imaging in medicine II: modelling and reconstruction, Physics in Medicine and Biology \textbf{42} (\textbf{5}) (1997), 841.

\bibitem{ArL} S. R. Arridge, W.R.B Lionheart, Nonuniqueness in diffusion-based optical tomography, Optics Letters \textbf{23}  (\textbf{11}) (1998), 882-884.

\bibitem{ArSc} S. R. Arridge and J.C. Schotland, Optical tomography: forward and inverse problems, Inverse Problems \textbf{25} (\textbf{12}) (2009), 123010.



\bibitem{Be-Fr} E. Beretta and E. Francini, Lipschitz stability for the electrical impedance tomography problem: the complex case, Communications in Partial Differential Equations \textbf{36} (2011), 1723-1749.

\bibitem{Be-dH-Fr-V-Z} E. Beretta, M. V. de Hoop, E. Francini, S. Vessella and J. Zhai, Uniqueness and Lipschitz stability of an inverse boundary value problem for time-harmonic elastic waves, Inverse Problems \textbf{33} (\textbf{3}) (2017), 035013.


\bibitem{C} A. P. Calder\'{o}n, On an inverse boundary value problem, Seminar on Numerical Analysis and its Applications to Continuum Physics (Rio de Janeiro, 1980),   65--73, Soc. Brasil. Mat., Rio de Janeiro, 1980. Reprinted in: Comput. Appl. Math. \textbf{25}  (\textbf{2-3}) (2006), 133-138.


\bibitem{DaHMa} B. Davey, J. Hill and S. Mayboroda, Fundamental matrices and Green matrices for non-homogeneous elliptic systems, Publicacions Matematiques \textbf{62} (\textbf{2}) (2018), 537 - 614.

\bibitem{dHQS} M. V. de Hoop, L. Qiu and O. Scherzer, Local analysis of inverse problems: H\"older stability and iterative regularization, Inverse Problems \textbf{23} (2012) 045001 (16pp).

\bibitem{DM} G. Dolzmann and S. M\"uller, Estimates for Green's matrices of elliptic systems by $L^p$ theory, Manuscripta mathematica \textbf{88} (1995), 261 - 273.




\bibitem{F} M. Fuchs, The Green-matrix for elliptic systems which satisfy the Legendre-Hadamard condition, Manuscripta mathematica \textbf{46} (1984), 97 - 115.

\bibitem{F1} M. Fuchs, The Green matrix for strongly elliptic systems of second order with continuous coefficients, Zeitschrift f\"ur analysis und ihre anwendungen \textbf{5} (6) (1986), 507 - 531 .





\bibitem{G} R. Gaburro, Stable determination at the boundary of the optical properties of a medium: the static case, Rend. Istit. Mat. Univ. Trieste \textbf{48} (2016), 407 - 431.

\bibitem{GL} R. Gaburro and W. R. B. Lionheart, Recovering Riemannian metrics in monotone families from boundary data, Inverse Problems \textbf{25} (\textbf{4}) (2009), 045004 (14pp).

\bibitem{G-S} R. Gaburro and E. Sincich, Lipschitz stability for the inverse conductivity problem for a conformal class of anisotropic conductivities, Inverse Problems \textbf{31} 015008 (2015).



\bibitem{Ha} B. Harrach, On uniqueness in diffuse optical tomography, Inverse Problems \textbf{25} (\textbf{5}) (2009), 055010.

\bibitem{HS} J. Heino and E. Somersalo, Estimation of optical absorption in anisotropic background, Inverse Problems \textbf{18} (\textbf{3}) (2002), 559-573.

\bibitem{HAS} J. Heino, S. Arridge, J. Sikora, and E. Somersalo, Anisotropic effects in highly scattering media, Physical Review E \textbf{68} (\textbf{3}) (2003), 031908.



\bibitem{H} N. Hyv\"onen, Characterizing inclusions in optical tomography, Inverse Problems \textbf{20} (\textbf{3}) (2004), 737.

\bibitem{KVKaAr} V. Kolehmainen, M. Vauhkonen, J. P. Kaipio and S. R. Arridge, Recovery of piecewise constant coefficients in optical diffusion tomography, Optic Express \textbf{7} (13) (2000), 468 - 480.










\bibitem{I1} V. Isakov, On the uniqueness in the inverse conductivity problem with local data, Inverse Problems and Imaging \textbf{1} (\textbf{1})  (2007), 95 - 105.













\bibitem{Mo} C. B. Morrey, Multiple integrals in the calculus of variations, Springer, Berlin (1966).







\bibitem{Sa} M. Salo, Inverse problems for nonsmooth first order perturbations of the Laplacian, PhD Thesis, University of Helsinki, Helsinki (2004).







\end{thebibliography}
\end{document}